# Beyond recency and magnitude: Learning-aware uncertainty estimation for safety stock under evolving forecast-error distributions

**Authors:** Luis Fernández-Palacios [a,*], Manuel Ceballos [b], Yolanda Muñoz-Ocaña [c]

[a] *Loyola University Andalusia, Department of Business Management, Campus Sevilla: Avda. de las Universidades, 2, 41704 - Dos Hermanas, Sevilla, Spain. E-mail address: lfpalacios@uloyola.es*

[b] *Loyola University Andalusia, Department of Engineering, Campus Sevilla: Avda. de las Universidades, 2, 41704 - Dos Hermanas, Sevilla, Spain. E-mail address: mceballos@uloyola.es*

[c] *Loyola University Andalusia, Department of Business Management, Campus Córdoba: Calle Escritor Aguayo, 4, 14004 Córdoba, Spain. E-mail address: yolandam@uloyola.es*

[*] *Corresponding author. Calle Estrellita Castro, 5, 41002 – Sevilla, Spain. E-mail address: lfpalacios@uloyola.es*

# Abstract

We develop a non-parametric approach to refine forecast-error histories for safety stock estimation without relying on error magnitude or recency alone. Using forecast-error data from a fast-moving consumer goods environment, we show that one year of errors is insufficient to characterize service-level risk, while pooling several years is problematic because error distributions change over time. Non-parametric annual comparisons indicate a progressive reduction in bias and dispersion, consistent with learning in the forecasting process. We propose LOWDII (Leave-One-Out Wasserstein Distributional Influence Index), a diagnostic that evaluates the distributional influence of each historical forecast error on the empirical uncertainty distribution used for safety stock calibration. LOWDII identifies observations whose influence is disproportionate to their representativeness, separating transient distortions from persistent tail behaviour without imposing parametric assumptions. We evaluate the method in a discrete-event simulation using subsequent-year demand realizations as validation. The results show that LOWDII achieves the target service level while reducing average stock by 5.1% to 22.6% relative to benchmark methods. From a managerial perspective, LOWDII helps firms translate improvements in the forecasting process into safety stock decisions, avoiding the projection of obsolete historical errors into the future and capturing the economic benefit of lower uncertainty requirements earlier.

# 1. Introduction

Safety stock remains economically relevant for fast-moving consumer goods (FMCG) companies because it serves as a buffer against operational volatility. The main advantage of using forecasts in production planning lies in improving the temporal positioning of supply decisions. Unfortunately, "All that we can say for certain about a forecast of demand is that it will be in error" (Silver et al., 1998, p. 109). Feng et al. (2011) argue that point demand forecasts only represent expected demand and that additional information on demand uncertainty, should be incorporated into safety stock calculations. This view is consistent with standard supply-chain formulations, where safety inventory is linked directly to forecast-error variability rather than to demand variability alone (Chopra & Meindl, 2016 p. 344). Saoud et al. (2025) argue that forecast uncertainty is central to inventory-related costs because forecast errors propagate into upstream supply-chain decisions. The pandemic altered purchasing patterns so profoundly that naïve forecast errors increased significantly (E2open, 2024), prompting firms to adopt inventory practices based on larger buffer stocks, thereby tying up substantial working capital (Song et al., 2025). In the company examined in this study, supply lead time variability represented only 3.4% and 3.8% of total variance in the historical analyses conducted. This is consistent with mature industrial settings, where firms often operate in complex and volatile markets shaped by aggressive commercial dynamics. Likewise, Barros et al. (2021) found that only 10% of the reviewed studies addressed supply uncertainty exclusively, whereas 30% examined it jointly with demand uncertainty.

Yet the quality of any safety stock estimate depends not only on the formula used, but also on whether the historical error sample is truly representative of the variability that the system should protect against. Arora et al. (2025) show that traditional forecasting methods may become unreliable during and after major disruptions because they typically rely on historical patterns whose continuity can be broken by shocks such as pandemics, geopolitical events, or sudden shifts in consumer behaviour. The same logic applies to safety stock estimation: historical forecast errors are useful only insofar as they represent future uncertainty. This representativeness problem involves a tension between statistical sufficiency and historical relevance. A short error history may provide too few replenishment cycles to estimate high service-level risk reliably. Conversely, a longer history may include errors generated under earlier forecasting regimes, temporary commercial shocks, or non-recurrent promotional conditions.

In the FMCG dataset analysed in this study, this tension is empirically visible. The company's current one-year forecast-error history is statistically fragile for high-service-level calibration, while direct multi-year pooling is also problematic because annual forecast-error distributions are not homogeneous. Non-parametric comparisons indicate differences across years in location and distributional shape, and longitudinal evidence suggests that many SKUs experience progressive centring of errors and reduction of dispersion over time. This pattern is consistent with learning in the forecasting process, but it also means that some historical errors may no longer be representative of the uncertainty structure that should inform future inventory protection.

To address this problem, this article proposes LOWDII (Leave-One-Out Wasserstein Distributional Influence Index), a non-parametric framework for assessing the representativeness of forecast-error histories before they are used in safety stock estimation. LOWDII treats each historical error as a candidate contributor to the empirical uncertainty distribution and evaluates whether its inclusion improves or distorts the distributional basis on which inventory protection is calculated. Its objective is not to minimize variability, but to retain variability that is structurally informative while excluding observations whose influence is disproportionate to their relevance for the current forecasting process. For this reason, it is especially appropriate for forecasting processes that learn from past errors, instead of assuming that the historical record is fully valid or evaluating deviations solely by their magnitude, recency, or distance from the median.

## 1.1 Literature review

In inventory management, the reliability of safety stock decisions depends critically on how demand uncertainty and forecast-error uncertainty are represented. Early work stresses that applying safety stock theory under an incorrect assumption about the demand distribution can have severe consequences, making the identification of the appropriate distributional shape central to the determination of safety factors (Brown, 1977). The problem is compounded when only limited representative samples are available, since smaller samples require corrective adjustments that translate into higher safety stock requirements (Moeeni et al., 2012). More broadly, safety stock estimation is not only a statistical forecasting problem, but also an inventory-parameterization problem, because the representation of uncertainty has direct operational consequences (Syntetos et al., 2016). Inventory control performance can deteriorate when demand parameters are unknown and must be estimated from historical observations, highlighting that safety-stock parameterization depends on the quality and reliability of the information used to estimate uncertainty (Strijbosch & Moors, 2005). This concern is reinforced in dynamic operating environments, where demand volatility and changing customer needs can make previous stock deployment and supply-policy configurations inefficient, requiring periodic review of the historical information used to calibrate inventory parameters (Cantini et al., 2024).

The literature is well developed in modelling forecast-error distributions for safety stock estimation, but provides much less explicit guidance on how to construct a representative historical sample of forecast errors. In practice, smoothed versions of forecast-error measures are often used to improve responsiveness (Boylan & Syntetos, 2008), while other work evaluates the error distributions generated by different base and combined forecasting methods rather than relying on a single forecast-error sample (Barrow & Kourentzes, 2016). These contributions support the broader premise that safety stock performance depends on the consistency between the forecast-error distribution learned from historical data and the error distribution observed out of sample. Model-based approaches derive inventory levels from simulated predictive distributions instead of relying only on approximate standard-deviation-based calculations (Snyder et al., 2002), while volatility-based approaches model time-varying demand variance and show that temporal heteroscedasticity can materially affect inventory performance (Zhang, 2007). When normality and iid assumptions are not adequate, KDE- and GARCH-based methods have been used to estimate lead-time forecast-error uncertainty, and forecast-error distributions have also been modelled ex post by combining alternative quantile estimates from the same historical error sample (Trapero et al., 2019a; Trapero et al., 2019b). This empirical line has been extended by incorporating EVT to model the tails of forecast-error distributions generated by extreme events (Derbel et al., 2022), and by using bootstrap-based approaches to set safety stock levels under skewed, multimodal, or highly variable lead-time demand (Saldanha et al., 2023).

However, the methodological logic of most of these contributions operates after accepting the historical error sample as a valid input. The present article intervenes at an earlier stage. In this study, we show that forecast-error distributions exhibit substantial year-to-year variation in both bias and distributional shape as the forecasting system learns and improves. As a result, we seek to cleanse the sample of isolated observations with abnormal distributional influence while preserving structural variability. This complements the existing literature by focusing not on the analytical propagation or modelling of forecast-error variance, but on the representativeness of the empirical forecast-error history from which the uncertainty distribution is learned.

Outside inventory management, the construction and assessment of representative historical samples has been mainly addressed in applied statistical modelling, forecasting, statistical learning, and time-series analysis through methods that adapt either the amount or the weight of past data under non-stationarity. Adaptive-window approaches expand or shrink the retained history as concept drift is detected (Bifet & Gavaldà, 2007), while forecasting studies select historical windows using out-of-sample predictive performance to avoid estimating current models from obsolete regimes under structural instability or time-varying parameters (Hirano & Wright, 2022; Huang & Wang, 2025; Inoue et al., 2017). Related approaches replace hard window selection with temporal weighting and online distribution estimation (Stephanou et al., 2017), or address non-representative episodes through model-based treatment of extremes and crisis effects in time series (McElroy, 2016, 2026).

These approaches share the premise that historical observations should not be treated as homogeneous or equally representative when the underlying data-generating process changes over time. However, they mainly operate through window selection, temporal weighting, drift adaptation, or model-based adjustment of extreme periods. They rarely address whether a specific individual observation remains distributionally representative of the empirical uncertainty structure used for a downstream decision. This limitation motivates a closer connection with three related methodological traditions: influence diagnostics in time series, non-parametric distributional comparison, and robust statistics.

A first relevant stream concerns the identification of influential observations and outliers in time series. Early contributions adapted influence diagnostics to serially dependent data by evaluating how individual observations affect model parameters, fitted values, or subsequent forecasts (Peña, 1990; Lee & Hui, 1993). This line of research was later extended through procedures explicitly designed for outlier-contaminated forecasting models (Chen & Liu, 1993), local influence diagnostics for ARIMA-type models (Lu et al., 2012), and forecast-oriented criteria that assess the effect of atypical observations on fitted and predicted values rather than relying only on standardized residuals or classical regression diagnostics (Oesterreich, 2020). More recently, the impact of recurring outliers in time-series data has been characterized through influence functionals, showing that some observations may be particularly harmful for parameter estimation and future prediction independently of their frequency (Li et al., 2021). These studies provide an important precedent for treating observation-level influence as a relevant statistical object. However, their main target is usually the sensitivity of model parameters, fitted values, or forecasts, rather than the distortion induced by each observation on the empirical distribution of forecast errors itself.

A second relevant stream comes from non-parametric robustness and distributional comparison. Statistical depth provides a way to assess centrality and outlyingness without imposing a specific parametric distribution, and has been used to design robust outlier identifiers with explicit attention to masking and swamping effects (Zuo & Serfling, 2000; Dang & Serfling, 2010). In parallel, robust statistics has long emphasized the use of bounded-influence reasoning and high-breakdown location and scale measures, including the median and the median absolute deviation, when observations may be contaminated or heavy-tailed (Hampel et al., 1986; Rousseeuw & Leroy, 1987; Iglewicz & Hoaglin, 1993). Wasserstein distances add a complementary distributional perspective because they compare probability distributions through the amount of mass displacement required to transform one distribution into another, making them sensitive not only to location but also to dispersion, asymmetry, and tail structure (Villani, 2009; Panaretos & Zemel, 2019; Peyré & Cuturi, 2019). In a related non-parametric setting, optimal-transport-based distances provide a rigorous framework for comparing empirical distributions without reducing distributional differences to a single location or scale parameter (Hallin et al., 2021). This is particularly relevant when forecast-error samples may contain heavy tails, asymmetries, or local distributional features that are operationally meaningful for safety stock calibration. In forecasting contexts, Wasserstein-based measures have been used to quantify disagreement among probabilistic and density forecasts, thereby supporting the interpretation of forecast distributions as objects that can be compared directly in a metric space (Cumings-Menon et al., 2021). Nevertheless, this literature generally compares distributions across models, forecasters, or samples, while influence diagnostics usually focus on observation-level effects on model outputs.

## 1.2 Empirical evidence on historical-span insufficiency and distributional non-homogeneity.

Cycle service level measures the percentage of replenishment cycles in which no stockout occurs. If the target service level is 98%, the implied stockout probability is 2% per cycle. Therefore, the statistical precision of the estimated service level depends on the number of cycles available in the historical sample of forecast errors. That number of cycles can be calculated using the formula for the confidence interval of a proportion:

$$\hat{p} \pm z \cdot \sqrt{\frac{\hat{p}(1-\hat{p})}{n}}$$

Solving for $n$, we obtain:

$$n = \frac{z^2 \cdot \hat{p}(1-\hat{p})}{e^2}$$

A larger sample is needed as the required precision and reliability of the service-level estimate increase. Walpole et al. (2012, Chapter 9) present examples with e = 0.02 for quality studies, service level, or compliance percentages. And they use Z=1.96 as the standard critical value for a 95% confidence level in formal studies or papers. With a service level of 98%, the formula gives 188 cycles which, depending on lot size, would be equivalent to between 3 and 4 years. That is well above the current standard of one year.

We have six years of forecasts available to work with. Within each year we have the weekly forecast for the following 58 weeks per SKU. For Hyndman and Kostenko (2007), calendar years are relevant in a seasonal business; the span used should cover multiples of calendar years to include all seasonal events. We exclude the first year, 2020, from the analysis because it coincides with the COVID pandemic year, when market dynamics were dramatically altered. Chatfield (2000, Chapter 6) warns us that historical spans including structural breaks can bias model parameters and it is often preferable to restrict the span to post-change data. According to the E2open (2024) Benchmark, at the peak in the Spring of 2020, almost half of all volume was exposed to extreme error. We reserve the last year, 2025, to validate our SS proposal in a DES simulation, leaving 4 years to develop the statistics of our method—sufficient, as we saw.

Syntetos and Boylan (2006) argue that forecasting methods should not be evaluated solely on statistical accuracy, but also according to their impact on inventory performance. In the same line, Strijbosch et al. (2011) show that forecasting and stock-control objectives do not necessarily coincide. Finally, Barrow and Kourentzes (2016) emphasize that forecast evaluation for inventory management should not be limited to accuracy metrics, since the full shape of the forecast-error distribution directly affects safety stock calculation. This argument is reinforced by Saoud et al. (2025), who propose the Ratio of Forecast Uncertainty to analyse the bullwhip effect from the perspective of forecast-error behaviour rather than demand variability alone. This supports the view that forecast evaluation for inventory decisions should not rely only on average accuracy indicators, but should also consider the uncertainty structure generated by forecast errors.

In this sense, year-to-year homogeneity should be assessed not only in terms of accuracy, but also in terms of distributional location and shape. If consecutive years differ substantially in these dimensions, they cannot be assumed to provide equally representative evidence about the current forecast-error process. We therefore used the Kruskal–Wallis (KW) test to assess whether annual error distributions differ in location, with particular attention to shifts in the median. The median is less sensitive to extreme observations than the mean and provides an interpretable indication of whether different years exhibit comparable levels of optimism or pessimism in the forecast errors. We complemented this with the Kolmogorov–Smirnov (KS) test, which compares the empirical distributions more broadly and can detect differences in shape, dispersion, and tail behaviour. This is relevant for safety stock estimation because heavier or more dispersed tails may translate into higher protection levels.

The results obtained from these tests were not what we expected. Barely 15% of the SKUs analysed had homogeneous results in both tests for three or more years (taking 2024 as the reference) when we needed between 3 and 4 years. The KS test rejected distributional homogeneity in 44% of the comparisons, while the KW test showed even stronger differences: comparable error levels were found in only 25% of cases. The most surprising aspect is that high-rotation formats, with theoretically stable commercial behaviour, did not have more than 1 or 2 homogeneous years. Since one year is statistically insufficient and four years are distributionally heterogeneous, we need to identify which individual historical errors remain representative of the uncertainty structure that should inform safety stock calibration.

The medians and quantiles of these formats tended to move synchronously across years, converging towards the centre and reducing their amplitude over time. This suggests a damped learning process: if one year is excessively optimistic, the next tends to be less so, and heavy tails tend to narrow in subsequent years. In other words, forecast errors become progressively more centred, while dispersion declines. This interpretation is consistent with Beutel and Minner (2012), who show that safety stock can be parameterized from forecast errors generated by models incorporating external explanatory variables. Since such models account for observable demand drivers, they also support the view that forecast-error distributions may improve as forecasting processes become more effective. In practice, this improvement may arise not only from enriched statistical models, but also from organizational learning about promotions, commercial actions, and forecasting bias.

To assess longitudinal evidence of learning, we computed two robust annual indicators for each SKU: the absolute annual median forecast error, capturing central bias, and the annual interquartile range, capturing dispersion. Negative slopes in these indicators, estimated against a normalized time index, indicate progressive centring of the forecast-error distribution and narrowing of its central mass, respectively. Given the limited number of annual observations per SKU, slope signs were interpreted as directional evidence rather than formal significance tests. Approximately 90% of SKUs exhibited at least one of these behaviours, while 61% exhibited both, indicating that their forecast errors became both more centred and less dispersed over time. Rather than a strictly monotonic convergence process, the results reveal a downward tendency with fluctuations. Forecast errors tend to become more centred and less dispersed over time, but this improvement is not smooth: individual SKUs may experience temporary reversals, renewed dispersion, or changes in the dominant sign of the error due to promotional intensity, commercial strategy, market-share initiatives, or other evolving business conditions. Therefore, the empirical pattern is better described as damped oscillatory behaviour than as a linear decay of forecast uncertainty.

This dynamic has direct implications for how historical forecast errors should be used in safety stock estimation. If improvement were monotonic, a purely temporal weighting scheme could be sufficient: older errors would simply receive lower weight and recent errors higher weight. However, when learning is oscillatory and affected by temporary commercial or operational shocks, recency alone is not a reliable proxy for representativeness. A recent error may reflect a non-recurrent distortion, while an older error may still belong to a persistent tail structure that should remain part of the uncertainty basis. The relevant task is therefore not merely to identify which part of the past is most representative of the present, but to shape the empirical distribution that will be projected into the future. In this sense, the objective is to preserve structural variability while reducing the influence of observations that mainly transmit transient distortions from previous or non-recurrent forecasting regimes. Including them increases safety stock and may reduce the incentive to improve forecasts, because even extreme forecasting failures are absorbed by inventory.

## 1.3 Contribution

This study addresses a prior and underexplored problem in safety stock calibration: whether the historical forecast errors used to estimate uncertainty are representative of the future variability that the inventory system should protect against. The empirical diagnosis shows that a one-year forecast-error history is statistically fragile for estimating high service-level risk, whereas direct multi-year pooling is also unreliable because annual forecast-error distributions differ in location, dispersion, and shape. The problem solved here is therefore not how to improve the point forecast, nor how to impose a parametric model on forecast errors, but how to construct a representative empirical uncertainty distribution before safety stock is calculated.

To the best of our knowledge, this is one of the first studies to explicitly formulate the representativeness of historical forecast errors as a methodological problem in safety stock calibration, and to address it through an influence-based diagnostic criterion. Under this criterion, an observation is considered problematic when its removal produces an unusually large change in the empirical uncertainty distribution relative to the influence pattern of the remaining observations for the same SKU. LOWDII identifies the dependence of the empirical distribution on each observation, quantifies that dependence over the full distributional shape, and determines whether the resulting influence is unusually large within the SKU-specific error structure. It therefore does not treat non-representativeness as a simple function of recency, magnitude, or distance from the median, but as a distributional influence problem. This allows the method to retain recurrent tail behaviour while reducing the influence of transient errors generated under earlier or non-recurrent forecasting conditions.

Empirically, the study validates that one year of weekly errors is insufficient for stable high-service-level calibration, that four-year pooling is distributionally heterogeneous according to KW and KS tests, and that LOWDII achieves the target service level in subsequent-year discrete-event simulation while reducing average stock by 5.1% to 22.6% relative to benchmark methods. From a managerial perspective, LOWDII is readily applicable in FMCG settings where products have life cycles long enough to provide multi-year forecast-error histories and where the forecasting process is sufficiently stable and disciplined to generate comparable errors over time. In these conditions, the method allows firms to translate improvements in the forecasting process into safety stock decisions, avoiding the projection of obsolete historical errors into the future and capturing the economic benefit of lower uncertainty requirements earlier.

## 1.4 Major assumptions and modelling choices

The main assumptions and modelling choices of this study are summarized as follows. First, we focus on forecast-error variability as the primary source of uncertainty for safety stock calibration. This choice is consistent with the empirical setting analysed, where supply lead-time variability represented only a small share of total variance in the historical analyses conducted. Second, we assume that the available forecast-error history is sufficiently long to support distributional diagnosis, but not homogeneous enough to be pooled without assessment. This motivates the use of LOWDII as a representativeness diagnostic rather than as a conventional outlier-removal rule.

Third, the empirical evidence is interpreted as consistent with progressive but non-linear learning in the forecasting process. LOWDII does not model learning directly; rather, it operationalizes one observable implication of learning: some historical errors may cease to be representative of the uncertainty structure relevant for future safety stock calibration.

Fourth, LOWDII is developed as a non-parametric procedure. No assumption is made regarding normality, symmetry, unimodality, or homoscedasticity of forecast errors. Instead, each observation is evaluated according to its influence on the empirical forecast-error distribution.

Fifth, the system is modelled under a periodic-review, order-up-to replenishment logic in a two-echelon setting. The analysis does not explicitly model capacity constraints, risk pooling, product substitution, or multi-echelon inventory interactions. These choices make the experimental comparison more controlled, but also delimit the scope within which the results should be interpreted.

Finally, the empirical analysis uses complete calendar years as the relevant historical unit, since FMCG forecast errors may be shaped by annual seasonality, promotional calendars, and recurrent commercial events. The pandemic year is excluded as an exceptional structural disruption, and the final year is reserved for out-of-sample validation.

# 2. Proposed mathematical method: LOWDII

Kutner et al. (2004, Chapter 10) emphasize that atypical observations should not be assessed only by their extremeness, but by whether they are anomalous, influential, and capable of distorting inference. This distinction is central here. If the extreme value is a spurious observation generated by idiosyncratic circumstances, retaining it biases the estimation of variability upward and results in excessive safety stock by overextending the upper tail of the distribution. Conversely, if the value is a legitimate manifestation of recurrent but irregular promotional demand, removing it suppresses structurally relevant variation and leads to an underestimation of the inventory buffer required to absorb uncertainty shocks. That is why not only the distance of an outlier from the mode is relevant, but also how unusual that value is relative to the rest of the errors. This multimodality of the distribution has been reported by other authors such as Barrow and Kourentzes (2016) and Manary et al. (2009).

In the presence of heavy tails, structural asymmetries, regime changes, or multimodality—features that are common in forecast errors for FMCG products—conventional methods tend to over-detect legitimate points belonging to real tails (swamping), or to miss isolated but highly distorting outliers (masking). Robust literature has repeatedly shown these effects. We need a technique that estimates how much the distribution changes when we remove observation i, i.e., how unusual that value is relative to the rest. For that, we propose LOWDII a fully non-parametric approach based on leave-one-out distributional influence, which evaluates how much the empirical distribution changes when a single observation is removed. LOWDII is therefore not a conventional magnitude-based outlier-removal procedure.

Leave-one-out evaluation identifies the dependence of the empirical distribution on each observation, but does not define a distribution-sensitive measure of that dependence. Wasserstein distance provides such a measure by capturing how the entire empirical distribution, including tails and local shape, changes when an observation is removed, but it does not by itself establish whether the observed change is unusually large for that SKU. MAD-based standardization supplies this robust internal reference scale, but only after the relevant influence quantity has been defined. LOWDII is an influence-based diagnostic criterion: an observation is considered problematic when its removal produces an unusually large change in the empirical uncertainty distribution relative to the influence pattern of the remaining observations. In this context, LOWDII does not attempt to forecast organizational learning directly. Rather, it operationalizes one of its observable consequences: as forecast errors become progressively more centred and less dispersed, some historical observations cease to be representative of the uncertainty structure that should inform future safety stock.

A purely time-based method assumes that representativeness is mainly driven by recency, although the learning process may be trend-improving yet non-monotonic, with fluctuations, setbacks, and disruptions. A magnitude-based outlier rule assumes that non-representativeness is mainly driven by distance from the centre, although forecast errors may exhibit heavy-tailed behavior. LOWDII adopts a different criterion: it evaluates whether each observation remains distributionally coherent with the empirical error structure after learning has occurred. This allows the method to retain recurrent tail behaviour while reducing the influence of transient errors generated under earlier or non-recurrent forecasting conditions.

## 2.1 Empirical Distribution Function (EDF)

Let $s$ denote a given SKU and let

$$\mathcal{E}_s = \{e_{s,1}, e_{s,2}, \dots, e_{s,n_s}\}$$

be the historical sample of forecast errors available for that SKU over the calibration period. For notational simplicity, the SKU subscript is omitted in the remainder of this section, and the sample is written as

$$\mathcal{E} = \{e_1, e_2, \dots, e_n\}.$$

Each $e_i$ represents one historical forecast error observation. The objective is to determine whether each observation contributes structurally relevant information to the empirical uncertainty distribution, or whether its inclusion produces a disproportionate distortion in the distribution used for safety stock calibration.

The empirical distribution function associated with the forecast-error sample $\mathcal{E}$ is defined as

$$\hat{F}_n(x) = \frac{1}{n}\sum_{j=1}^{n} \mathbf{1}\{e_j \leq x\}$$

where $\mathbf{1}\{\cdot\}$ is the indicator function, equal to 1 if the condition inside the braces is true and 0 otherwise.

The EDF assigns probability mass to each observed forecast error and provides a non-parametric representation of the historical error distribution, without assuming normality, symmetry, unimodality, homoscedasticity, or any other specific distributional form (Van der Vaart, 1998). It therefore constitutes the empirical starting point from which the LOWDII influence index is constructed.

Let the ordered sample be

$$e_{(1)} \leq e_{(2)} \leq \cdots \leq e_{(n)}.$$

Then, at each ordered observation $e_{(k)}$, the EDF is

$$\hat{F}_n(e_{(k)}) = \frac{k}{n}, k = 1, \dots, n.$$

Equivalently, for any real value $x$,

$$\hat{F}_n(x) = \frac{|\{j: e_j \leq x\}|}{n}$$

This EDF is the empirical uncertainty distribution that would be used if the entire historical sample were retained. LOWDII evaluates whether each observation $e_i$ is distributionally coherent with this empirical distribution.

## 2.2 Leave-one-out empirical distribution

For each observation $e_i$, define the leave-one-out sample

$$\mathcal{E}^{(-i)} = \mathcal{E} \setminus \{e_i\} = \{e_1, \dots, e_{i-1}, e_{i+1}, \dots, e_n\}.$$

The empirical distribution function associated with this reduced sample is

$$\hat{F}_{n-1}^{(-i)}(x) = \frac{1}{n-1} \sum_{\substack{j=1 \\ j \neq i}}^{n} \mathbf{1}\{e_j \leq x\}.$$

The difference between $\hat{F}_n$ and $\hat{F}_{n-1}^{(-i)}$ captures how much the empirical distribution changes when observation $e_i$ is removed. If removing $e_i$ produces only a small change, then $e_i$ is not distributionally influential. If removing $e_i$ substantially changes the empirical distribution, then $e_i$ has high distributional influence.

## 2.3 Leave-one-out Wasserstein distributional influence

To quantify the change between the full EDF and the leave-one-out EDF, define the distributional influence of observation $e_i$ as

$$\Delta_i = D\left(\hat{F}_n, \hat{F}_{n-1}^{(-i)}\right)$$

where $D(\cdot,\cdot)$ is a distance between probability distributions.

In LOWDII, the selected distance is the first-order Wasserstein distance, denoted by $W_1$. Wasserstein distance, which measures "how much work" is needed to transform one distribution into the other. It captures changes across the whole distribution, including tails. It is sensitive to shape, not only to spikes.

Therefore,

$$\Delta_i = W_1\left(\hat{F}_n, \hat{F}_{n-1}^{(-i)}\right)$$

For a distribution function $F$, let $F^{-1}$ denote the generalized inverse, defined as

$$F^{-1}(u) = \inf\{x \in \mathbb{R}: F(x) \geq u\}, u \in (0,1).$$

For one-dimensional distributions, the first-order Wasserstein distance can be expressed as the integral of the absolute difference between their quantile functions:

$$W_1(F, G) = \int_0^1 | F^{-1}(u) - G^{-1}(u) | \, du.$$

Thus, for each observation $e_i$,

$$\Delta_i = \int_0^1 | \hat{F}_n^{-1}(u) - \left(\hat{F}_{n-1}^{(-i)}\right)^{-1}(u) | \, du$$

Equivalently, using the cumulative distribution representation in one dimension,

$$W_1(F, G) = \int_{-\infty}^{+\infty} | F(x) - G(x) | \, dx.$$

Therefore,

$$\Delta_i = \int_{-\infty}^{+\infty} | \hat{F}_n(x) - \hat{F}_{n-1}^{(-i)}(x) | \, dx$$

This formulation makes explicit that $\Delta_i$ measures the global distributional change induced by removing observation $e_i$. It is not a pointwise residual, nor a standardized distance from the mean or median. It is a full-distribution influence measure.

In practical terms:

$$\Delta_i \approx 0$$

means that removing $e_i$ barely changes the empirical forecast-error distribution. Conversely,

$$\Delta_i \gg 0$$

means that $e_i$ contributes disproportionately to the empirical distributional shape.

A central observation will generally produce a small $\Delta_i$, because removing it does not materially alter the tails or the global structure of the EDF. A recurrent tail observation may produce a moderate $\Delta_i$, because although it is extreme, it belongs to a region already supported by other observations. By contrast, an isolated observation located in a low-density region can produce a large $\Delta_i$, because it creates or extends a tail, generates an artificial local mode, or changes the apparent dispersion of the historical error distribution.

LOWDII does not penalize an observation simply because it is far from the centre. It penalizes an observation when its removal produces an unusually large change in the empirical distribution. This is consistent with the methodological objective stated in the manuscript: preserve structural variability while reducing the influence of transient or non-representative historical distortions.

## 2.4 Discrete empirical implementation of the raw influence measure

For computational implementation, it is useful to distinguish between the empirical distribution function and the empirical probability measure that generates it. The full empirical probability measure associated with the forecast-error sample $\mathcal{E}$ is

$$\hat{P}_n = \frac{1}{n}\sum_{j=1}^{n} \delta_{e_j},$$

where $\delta_{e_j}$ denotes a unit point mass located at observation $e_j$.

The corresponding empirical distribution function is

$$\hat{F}_n(x) = \hat{P}_n((-\infty, x]) = \frac{1}{n}\sum_{j=1}^{n} \mathbf{1}\{e_j \leq x\}.$$

Similarly, the leave-one-out empirical probability measure obtained after removing observation $e_i$ is

$$\hat{P}_{n-1}^{(-i)} = \frac{1}{n-1}\sum_{\substack{j=1 \\ j\neq i}}^{n} \delta_{e_j}.$$

Its associated empirical distribution function is

$$\hat{F}_{n-1}^{(-i)}(x) = \hat{P}_{n-1}^{(-i)}((-\infty, x]) = \frac{1}{n-1}\sum_{\substack{j=1 \\ j\neq i}}^{n} \mathbf{1}\{e_j \leq x\}.$$

The raw distributional influence of observation $e_i$ can therefore be written as the first-order Wasserstein distance between the two empirical probability measures:

$$\Delta_i = W_1\left(\hat{P}_n, \hat{P}_{n-1}^{(-i)}\right).$$

Equivalently, using the corresponding empirical distribution functions, the same quantity can be written as

$$\Delta_i = W_1\left(\hat{F}_n, \hat{F}_{n-1}^{(-i)}\right).$$

In the empirical implementation, the full and leave-one-out samples are treated as normalized unweighted empirical distributions, assigning masses $1/n$ and $1/(n-1)$, respectively. Hence, removing $e_i$ affects the empirical measure both by deleting its point mass and by renormalizing the remaining observations. The raw influence is computed as

$$\Delta_i = \text{wasserstein_distance}\left(\mathcal{E}, \mathcal{E}^{(-i)}\right), i = 1, \ldots, n.$$

This generates the raw influence vector

$$\boldsymbol{\Delta} = (\Delta_1, \Delta_2, \ldots, \Delta_n).$$

The resulting vector summarizes how influential each historical forecast error is for the empirical uncertainty distribution of the SKU.

## 2.5 Robust normalization of the raw influence values

The raw $\Delta_i$ values are SKU-specific. Their magnitude depends on the scale of the forecast errors, the dispersion of the SKU, the sample size, and the structure of the empirical distribution. Therefore, a direct comparison of raw $\Delta_i$ values across SKUs is not appropriate.

LOWDII standardizes the influence values internally, within each SKU, using a robust location and scale transformation.

Let

$$m_\Delta = \text{median}(\Delta_1, \ldots, \Delta_n)$$

be the median influence level for that SKU.

The median absolute deviation of the influence values is

$$MAD_{\Delta} = \text{median}(|\ \Delta_i - m_{\Delta}\ |: i = 1, \dots, n)$$

$$LOWDII_i = \begin{cases} \dfrac{\Delta_i - m_{\Delta}}{MAD_{\Delta}}, & MAD_{\Delta} > 0, \\ 0, & MAD_{\Delta} = 0. \end{cases}$$

Using MAD is consistent with the recommendations of Iglewicz and Hoaglin (1993) and with the robust tradition (Rousseeuw & Leroy, 1987), providing stability under heavy tails.

The robust Z-score for observation $e_i$ is then

$$LOWDII_i = \frac{\Delta_i - m_{\Delta}}{MAD_{\Delta}}$$

# 2.6 LOWDII decision rule

In the literature we have not found a consensus on recommended thresholds above which LOWDIIi allows excluding an observation. But the bibliography does describe a phenomenon: when a point is far away, its influence on a robust estimator increases; but when it is very far away, the transition from moderate to explosive influence is not smooth. Hampel et al. (1986) do not define a threshold but warn that 3.5 detects too many false positives and indicates that in asymmetric or heavy-tailed distributions the threshold should be increased. In robust detection frameworks (Python PyOD, R robustbase, ESD-type algorithms), "extreme anomalies" are usually between 6 and 8 MADs, and "impossible" beyond 8 MADs.

For this reason, the threshold was not selected solely from the literature, but also checked against the empirical behaviour of LOWDII in the dataset. Figure 1 ranks the LOWDII values within each SKU from highest to lowest. The x-axis represents the observation rank inside each SKU, while the y-axis reports the corresponding LOWDII value, calculated from the robustly standardized Wasserstein influence score. The plot shows that very high positive LOWDII values are concentrated in the first ranks, whereas the remaining observations quickly converge towards moderate, near-zero, or negative values. This indicates that distributionally influential observations are sparse and concentrated in a limited subset of records.

The sharp decay observed in the ranked curves supports the use of a high exclusion threshold. In particular, (LOWDII_i > 8) lies in the region where influence becomes clearly exceptional relative to the SKU-specific distribution of LOWDII values. The exclusion rule is therefore one-sided: only observations with (LOWDII_i > 8) are removed. Negative values are retained because they represent lower-than-median distributional influence, not excessive distortion.

**Figure 1. Ranked LOWDII influence scores by SKU**

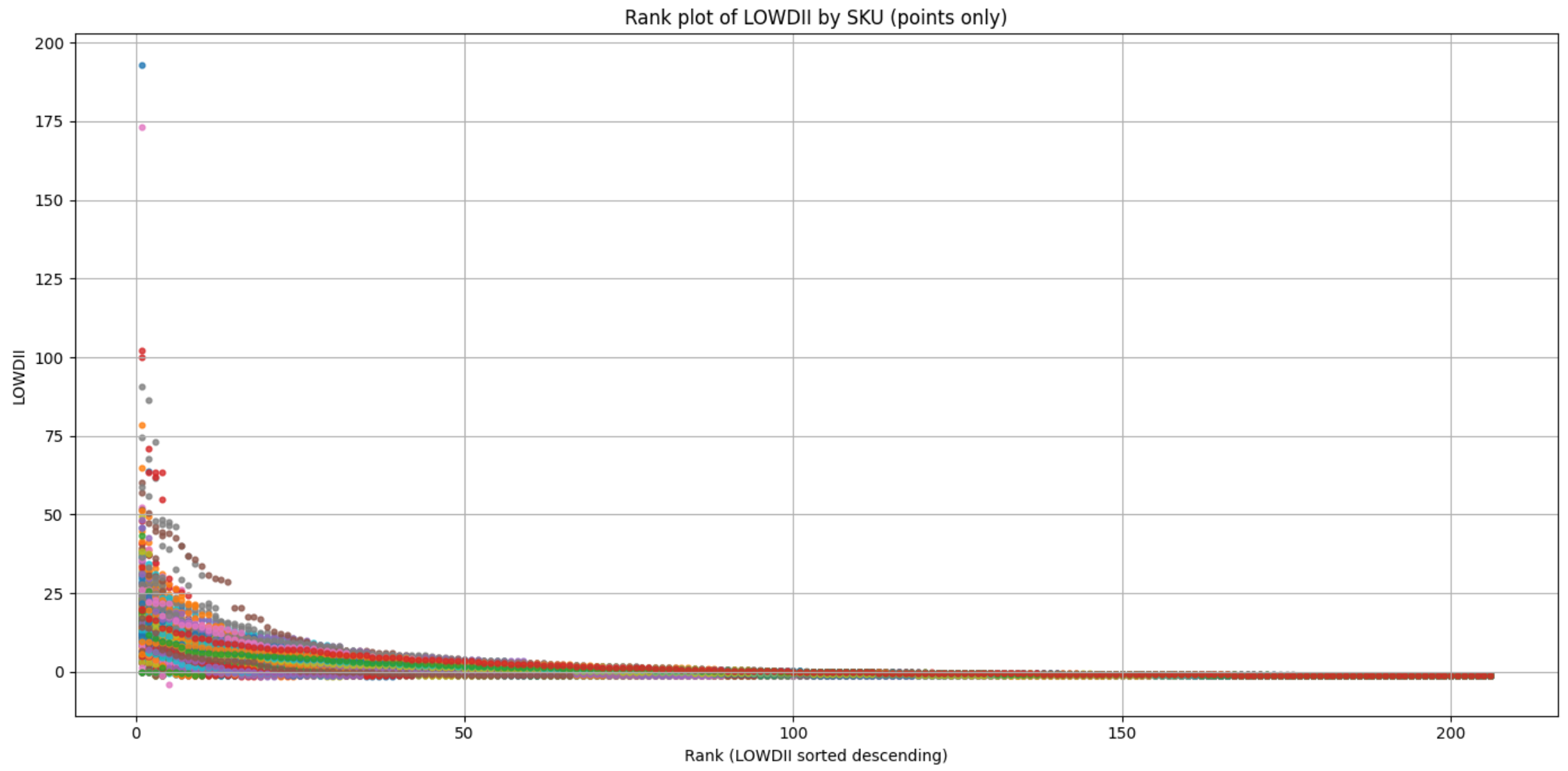

Let $\tau$ denote the influence threshold.

$$\tau = 8,$$

based on the observed explosive behaviour of ranked $LOWDII_i$values around that reference level.

Observation $e_i$ is classified as distributionally incompatible when

$$LOWDII_i > 8.$$

Therefore, the LOWDII exclusion indicator is

$$R_i = \mathbf{1}\{LOWDII_i > \tau\}.$$

The cleaned historical error sample is then

$$\mathcal{E}^{LOWDII} = \{e_i \in \mathcal{E}: LOWDII_i \leq \tau\}.$$

The empirical distribution used for safety stock calibration becomes

$$\hat{F}^{LOWDII}(x) = \frac{1}{n^{LOWDII}} \sum_{e_i \in \mathcal{E}^{LOWDII}} \mathbf{1}\{ e_i \leq x \},$$

where

$$n^{LOWDII} = \sum_{i=1}^{n} \mathbf{1}\{ LOWDII_i \leq \tau \}$$

is the number of observations retained after LOWDII cleansing.

## 2.7 Safety stock calibration after LOWDII cleansing

Once the LOWDII-cleaned sample has been obtained, the safety stock calculation is based on the empirical variability of the retained forecast errors.

Let

$$\sigma_{LT}^{LOWDII}$$

denote the standard deviation of the cumulative lead-time forecast error estimated from $\mathcal{E}^{LOWDII}$. Then the safety stock is calculated as

$$SS^{LOWDII} = z_\alpha \cdot \sigma_{LT}^{LOWDII},$$

where $z_\alpha$ is the standard normal service factor associated with the target service level $\alpha$. For a 98% cycle service target,

$$\alpha = 0.98, z_\alpha \approx 2.054.$$

The standard deviation-based safety stock expression is retained deliberately because the experimental scope of this study is intentionally limited in order to isolate the incremental value of the LOWDII procedure. The objective is not to replace the conventional inventory formula, but to improve the representativeness of the empirical forecast-error sample from which the uncertainty parameter is estimated. This choice also ensures comparability with the company's current calibration logic and with the benchmark methods evaluated in the simulation. LOWDII therefore acts as a preprocessing procedure: it modifies the historical error basis used to estimate cumulative lead-time forecast-error variability, while leaving the downstream safety stock parameterization unchanged.

# 3. The impact of distributional representativeness on safety stock performance: simulation experiments

## 3.1 Simulation design and benchmark methods

The simulation compares LOWDII against five benchmark methods using the same forecast-error database, experimental units, performance metrics, and statistical testing framework

**Database and forecast-error construction.** The six methods will use the same database of weekly forecast error. According to Boylan and Babai (2016), lead-time forecasts can be interpreted as a form of overlapping temporal aggregation of demand. Prak et al. (2017) show that conventional procedures based on multiplying one-period forecast-error variance by the lead time may be flawed, because future forecast errors can be correlated across the lead-time interval. The results of Babai et al. (2022) show that not considering the level of aggregation, particularly under positive demand autocorrelation, may lead to higher inventory costs and lower service levels. For that reason, historical errors have been transformed into the empirical variance of the cumulative error in the lead time (Lee, 2014).

**Calibration and validation windows.** As explained above, the period 2021–2024 is used to calibrate the competing methods, providing a sufficiently long sample to develop the statistics of LOWDII. The final year, 2025, is reserved for out-of-sample validation through the discrete-event simulation described in Section 3.2.

**Experimental unit.** The SKU is used as the experimental unit. The final simulation sample contains 153 SKUs.

**Benchmark methods.** Six methods are compared.

- **No-cleansing baseline**: RAW. 4 years with no cleansing applied
- **Recency-based**:
  - SPAN: Only 1 year with no cleansing applied. The current method of the company under study
  - SMOOTH52 & SMOOTH208: Exponential smoothing estimates SKU-level local volatility by recursively updating the mean and variance of forecast errors, giving more weight to recent observations according to a half-life parameter (h). A 52-week half-life makes the estimate more responsive to recent changes, while a 208-week half-life produces a more stable long-run volatility measure but is less representative of the current forecasting regime. Since no stable seasonality was found in the errors, the method uses exponential weighting rather than Holt or Holt–Winters seasonal specifications.

- **Magnitude-based**: IQR. 4 years with the IQR method, the classical outlier detection technique based on the interquartile range, defined as the difference between the third and first quartiles: $IQR = Q_3 - Q_1$. Observations are typically flagged as outliers when they fall below $Q_1 - 1.5 \cdot IQR$ or above $Q_3 + 1.5 \cdot IQR$, providing a simple and robust rule that does not rely on normality assumptions.
- **Distributional-influence-based**: LOWDII. 4 years of cleansing using a distributional influence-based approach, explained in Section 2

**Performance metrics.** The methods are evaluated using four simulation outputs: average weekly inventory (AvST), beta case fill rate (BFR), percentage of weeks without stock-out (PctWeeks), and average days out of stock (AvOutDays). These metrics are defined in Section 3.3.

**Statistical comparison.** Since all methods are evaluated on the same SKUs, performance comparisons follow a repeated-measures design. For each metric, we apply repeated-measures ANOVA, Friedman tests, Holm-corrected pairwise comparisons, and descriptive statistics; the rationale and interpretation of these tests are detailed in Section 3.3.

## 3.2 Discrete-event simulation model

The simulation is designed to test whether the distributional representativeness criterion introduced by LOWDII translates into inventory performance. Discrete-event simulation is used because the inventory consequences of alternative safety stock policies depend on the interaction between forecast-error uncertainty, planning parameters and operational replenishment logic. This modelling choice is consistent with prior simulation-based studies of production planning systems. Gansterer et al. (2014), for example, use a discrete-event simulation framework to analyse the joint effect of planned lead times, safety stocks and lot sizes in hierarchical production planning. Similarly, Seiringer et al. (2022) examine how forecast error and forecast bias affect safety stock in an MRP system with rolling-horizon forecast updates. In the present study, the simulation is not used to optimise all planning parameters simultaneously, but to reproduce as faithfully as possible the operational planning logic of the company from which the data were obtained, allowing the inventory consequences of alternative historical forecast-error treatments to be evaluated under realistic batch scheduling logic, rolling-horizon forecast updates and service-level conditions.

As described above, 2021–2024 are used for calibration and 2025 for out-of-sample DES validation. Our DES is characterized by:

3.2.1 **Multiperiod Case:** This is a dynamic inventory model in which the system is analysed over multiple time periods using weekly buckets, allowing the evolution of inventory, demand and replenishment decisions to be captured over time. Each week, the model plans production for the week immediately after the production lead time (week t+5), prioritizing those SKUs whose gap between projected stock and safety stock is the smallest. Planned quantities comply with the minimum lot size that optimizes the current scheduling system.

3.2.2 **Projected stock:** Projected stock for week t+5 is calculated using the initial stock of the current week, minus the forecast for the next 5 weeks (the real forecast that existed in that week's plan), minus the production of the 4 weeks in the frozen period.

3.2.3 **Initial stock:** For each period, opening stock equals the previous week's stock plus production minus the previous week's actual demand.

3.2.4 **Continuous review approximation:** The simulation uses an order-up-to (OUT) logic to generate planned production orders aimed at replenishing finished-goods inventory over the lead-time planning horizon. The OUT policy is widely adopted in industry and is readily available in many ERP/MRP systems (Disney et al., 2025). In real operations, the production plan is reviewed weekly, typically on Mondays, which corresponds to a periodic-review scheme. However, the DES works with weekly buckets and collapses each week into a single event in which production arrivals occur before demand. In practice, this formulation approximates a continuous-review system, since it does not explicitly represent the *undershoot* effect associated with intra-week variability or with variability in the production day. This does not invalidate the model, because its objective is to measure the efficiency of *safety stock* to buffer the variability derived from forecast error. Nevertheless, later it will be necessary to incorporate an adjustment factor that captures that *undershoot* effect. In a classic periodic-review framework in continuous time (or a discretization where demand occurs throughout the interval), inventory must cover until the next decision can materialize; we should follow an approximation such as the one proposed by Campbell (1995) where our Protection Interval would be equivalent to the lead time (4) + the order interval (1), i.e., the planning frequency. We could also use a more sophisticated approximation depending on the probability of producing earlier or later within the week.

3.2.5 **WARM UP period:** The first week of the simulation corresponds to week 50 of 2024 and includes a frozen production period designed to ensure that, by the beginning of week 02 of 2025, all SKUs are above their safety stock (SS) levels. These four weeks are excluded from the analysis. In the first week after the warm-up period, each SKU starts with at least its corresponding SS level, and each group of SKUs sharing the same production resource (production line) holds cycle stock equal to the sum of half the lot size of each SKU.

3.2.6 **52 chained replicas of a 1-year simulation with weekly buckets:** For low-rotation items, where the production lot may cover more than one week of demand, service-level distortions depend strongly on the point of the replenishment cycle at which the disruption occurs. If the disruption takes place immediately after a production event, its impact is more likely to be absorbed by the available stock. In addition, these items complete fewer replenishment cycles per year. To ensure that disruptions and heavy-tail demand realizations are evaluated over all possible positions within the cycle, the model uses 52 chained replicas of a one-year simulation, interpreted as 52 consecutive years and aligned with the 52 weeks of the year.

3.2.7 **No lead time, no yield uncertainty in production:** As explained in Section 1., this effect represents less than 5% of total uncertainty and therefore is not a detailed object of this study.

3.2.8 **Single step simulation:** We only consider finished-goods inventory planning (not components). We could incorporate other levels by applying the same techniques developed in the proposal, but for model validation we simplify the test.

3.2.9 **Only produced SKUs:** We will work only with manufactured products, not imported items nor co-packings due to their lower representativeness. These represent only 2.4% of sales volume.

3.2.10 **No obsolete items:** Since we do not include innovations due to their short Span of History, we will not consider product obsolescence. This product typology is responsible for 90% of the special cases.

3.2.11 **No allocation dimensioning:** Although the company operates with multiple production and logistics sites, the simulation concentrates inventory decisions in a single location. This simplification is justified by the operational setting of the case study: replenishment lead times across sites are short, approximately 48 hours, and allocation decisions are therefore not the main source of uncertainty analysed in this research. Instead, the model focuses on the finished-goods inventory decision associated with forecast-error uncertainty and safety stock calibration. This is consistent with recent modelling practice where simulated supply chain is represented as a two-echelon system, as in Disney et al. (2025). This abstraction allows the analysis to isolate the effect of alternative historical error treatments without confounding it with multi-site allocation rules

3.2.12 **Batch scheduling and lot-size logic:** Lot-sizing is explicitly represented in the simulation because safety stock does not operate independently from production-planning parameters. Altendorfer (2019) shows that lot size, safety stock, planned lead time and capacity utilisation interact in stochastic multi-item production systems, affecting optimal safety stocks. This prevents the simulation from assuming an unrealistic continuous replenishment response and allows the effect of the safety stock policy to be evaluated under the actual batching constraints faced by planners.

3.2.13 **Controlled capacity setting:** In reality, supply crises may occasionally occur because of capacity shortages. However, correctly sizing the required production capacity in the medium and long term would constitute a different research line. To isolate the effect of the safety stock policy, the forecast-error treatment and the batch scheduling logic, weekly capacity for each group of SKUs sharing the same production resource is set equal to that group's real sales in the corresponding week. This modelling choice prevents stock variations from being driven by seasonal or structural capacity shortages, while still preserving the operational impact of demand volatility, safety stock levels and minimum lot-size constraints.

3.2.14 **Back ordering with no lost sales**: In our model, a stock-out occurs if, at the end of each period, the resulting stock after subtracting sales and adding production is negative. We work under a backordering-with-no-lost-sales system: the negative stock of each week accumulates and must be compensated in the following month if we want to avoid incurring greater service losses. We consider this assumption reasonable because the vast majority of customers keep stock in their own warehouses; in the business reality it has been observed that only when lack of availability extends beyond one week may there be a possibility of lost sales.

3.2.15 **No rationing game**: We exclude massive stock-outs (affecting all customers) caused by supply crises; we assume that the high number of customers will share stock-outs randomly in each period, generating a risk-pooling effect (Lee et al., 1997).

## 3.3 Performance metrics and statistical comparison tests

We measure the performance of our proposal through two primary indicators:

- **Average weekly inventory in €**: The objective is to minimize the average level.
- **Beta service level (BFR, %):** The beta service level, or fill rate, is the service indicator used by the company under study. It is computed as the ratio between units served and units requested. In this study, it is evaluated both at aggregate level and at SKU level to determine whether the simulated inventory policy achieves the target service level of 98%.

We must bear in mind that, in human relationships, there is not only the cold mathematics of results, but also the perception of them. We may meet the objectives stated above, but still have failed our customer too many times, or only once but for too long. These perceptions can later negatively affect the relationship with the customer and commercial negotiations. The system we propose must guarantee a BFR, but must not significantly harm the other two complementary metrics. Therefore, and although as secondary indicators, we include other metrics in our analysis.

- **Percentage of weeks without stock-out:** Blindly targeting BFR in high-rotation products can lead us to meet the indicator while experiencing continuous micro stock-outs.
- **Average duration of a stock-out in days:** It is not the same to fail 3 times in a year by small quantities as to fail only once but supply absolutely nothing for a prolonged period, especially if during that period some special action had been negotiated where both parties had invested money.

For each metric, we run:

1. Repeated-measures ANOVA
2. Friedman test
3. Pairwise post hoc comparisons across all method pairs
4. Holm correction for p-values
5. Descriptive statistics by method

We use this analytical framework because we compare multiple methods on the same SKUs, so observations are not independent and the design is repeated measures. For each metric, repeated-measures ANOVA tests whether there is an overall method effect while accounting for the paired structure of the data. As a robustness check, we also include the Friedman test, which does not require normality and assesses whether the conclusions hold from a non-parametric perspective. Once a significant overall effect is detected, we perform pairwise post hoc comparisons across all method pairs to identify which treatments truly differ from one another. Because this process multiplies the number of tests and increases the risk of false positives, we apply Holm's correction to the p-values, which controls the Type I error rate with more power than Bonferroni. Finally, we report descriptive statistics by method to complement inference with a substantive reading of magnitudes, dispersion, and effect direction.

## 3.4 Simulation results

At the aggregate level, considering the 52 replicas, the LOWDII method achieves the target service level of 98% in both BFR and PctWeeks while maintaining the lowest average stock. The largest differences are observed relative to RAW and SM208, the two approaches in which older observations retain greater weight. This supports the need for cleansing and is consistent with the system dynamics described through the Damped Oscillatory Behaviour. IQR also performs better than SM52, which may suggest that the relative position of observations is even more important than the time at which they were recorded. Finally, the SPAN method, which relies only on the most recent year of data, delivers results similar to those of LOWDII (Table 1). These results support the central claim that cleansing based on distributional influence can reduce inventory investment without sacrificing the target service level.

**Table 1. Aggregate level results compared to LOWDII**

| | LOWDII | Dif Vs RAW | Dif Vs SPAN | Dif Vs IQR | Dif Vs SM52 | Dif Vs SM208 |
|---|---|---|---|---|---|---|
| AvST | 333,034 | 19.0% | 0.2% | 5.1% | 17.5% | 22.6% |
| BFR | 98.1% | 0.9% | 0.2% | 0.3% | 0.8% | 0.1% |
| AvOutDays | 7.1 | -46.0% | -9.0% | -16.6% | -43.2% | -3.4% |
| PctWeeks | 97.8% | 1.2% | 0.1% | 0.3% | 1.0% | 0.5% |

# 3.5 Statistical significance of performance differences

At the SKU level, both the repeated-measures ANOVA and the Friedman test indicate significant overall differences across methods for all four performance metrics. This confirms that method choice has a systematic effect on inventory and service outcomes, under both parametric and non-parametric repeated-measures testing. The post hoc analysis with Holm correction shows that some method pairs differ significantly, whereas others do not. Taking LOWDII as the reference, the pattern is quite clear and can be observed in Table 2.

**Table 2. Post hoc with Holm correction results comparing LOWDII with the other 5 Methods**

| Metric | RAW | IQR | SMOOTH52 | SMOOTH208 | SPAN |
|---|---|---|---|---|---|
| **AvST** | **YES** (-412.9) | **YES** (-111.9) | **YES** (-381.3) | **YES** (-491.4) | **NO** |
| **BFR** | **YES** (-0.0060) | **YES** (-0.00189) | **YES** (-0.00535) | **NO** | **NO** |
| **AvOutDays** | **YES** (+0.0244) | **YES** (+0.0081) | **YES** (+0.0216) | **NO** | **NO** |
| **PctWeeks** | **YES** (-0.0114) | **YES** (-0.00305) | **YES** (-0.00963) | **YES** (-0.00458) | **NO** |

The results suggest that LOWDII improves safety stock calibration not by mechanically reducing variability, but by modifying the empirical uncertainty basis used in the calculation. Its advantage lies in distinguishing observations that remain structurally informative for future uncertainty from observations whose distributional influence is disproportionate to their representativeness. This allows inventory reductions to be achieved without a material deterioration in service performance, whereas the other methods, with the exception of SPAN, tend to achieve comparable service levels only at the cost of higher stock levels.

LOWDII and SPAN do not differ significantly in any of the four metrics. This suggests that, for these data, both methods exhibit very similar overall behaviour. Although SPAN shows performance statistically comparable to that of LOWDII across the evaluated metrics, its reliance on only one year of historical data makes it inherently less robust from a statistical perspective.

As shown for $e = 0.02$, $Z = 1.96$, and a target service level of 98%, parameter estimation based on a longer demand history is generally more reliable, with between three and four years of data typically required to achieve stable safety stock calibration. Therefore, the similarity between SPAN and LOWDII should be interpreted with caution.

## 3.6 Interpretation: representativeness versus recency, IQR, smoothing, and pooling

Because SPAN uses only the most recent year of data, it is more vulnerable to sampling variability and year-specific effects. By contrast, LOWDII achieves similar results using a longer historical horizon and a more robust cleansing procedure, making it a more stable and trustworthy approach for implementation.

The results also have relevant implications regarding how historical forecast errors should be interpreted before being used for safety stock calibration. Methods such as IQR act directly on the magnitude of the error. As a consequence, a single extreme error can shift the empirical quantile and lead the method to retain isolated observations that inflate the estimated uncertainty. LOWDII, by contrast, evaluates whether an observation is structurally consistent with the surrounding error distribution, allowing it to distinguish between relevant tail behaviour and observations that are exceptional because they are not in harmony with the rest of the historical pattern.

Time-based methods, such as SPAN and smoothing approaches, address the problem from a different angle by assigning greater importance to recent observations. However, forecast learning and error reduction are not necessarily linear over time. After several years of progressively smaller and more centred errors, new distortions may arise due to changes in commercial strategy, promotional intensity, or market-share capture initiatives. This can be observed in Figure 2, where the comparison between the in-sample period and 2025 for two specific SKUs shows that the out-of-sample distribution is not simply a continuation of the previous error pattern. In particular, the right-hand side of the superimposed histograms reflects the impact of a reactive and successful promotional action in 2025 that was not incorporated into the forecasting process, generating large positive forecast errors and altering the shape of the out-of-sample error distribution.

From an operational perspective, this is important because a method based only on recency may assign excessive weight to recent distortions, even when those observations correspond to conditions that are unlikely to recur in the same form or that the organization is expected to correct through forecast learning. LOWDII offers a more selective preprocessing logic: it can remove observations whose exceptional nature arises from their lack of coherence with the broader error structure, rather than from their mere distance from the central part of the distribution. In practical terms, these observations may correspond to errors generated under specific commercial or operational conditions that should not be automatically embedded into future safety stock calculations. By not assuming that all historical forecast errors remain representative of future uncertainty, LOWDII allows firms to capture the economic benefits of improved forecasting earlier, without compromising the target service level.

**Figure 2. Comparison between the in-sample period (2021-2024) and 2025 for two selected SKUs, illustrating the risk of assigning excessive weight to recent but potentially non-recurrent forecast errors.**

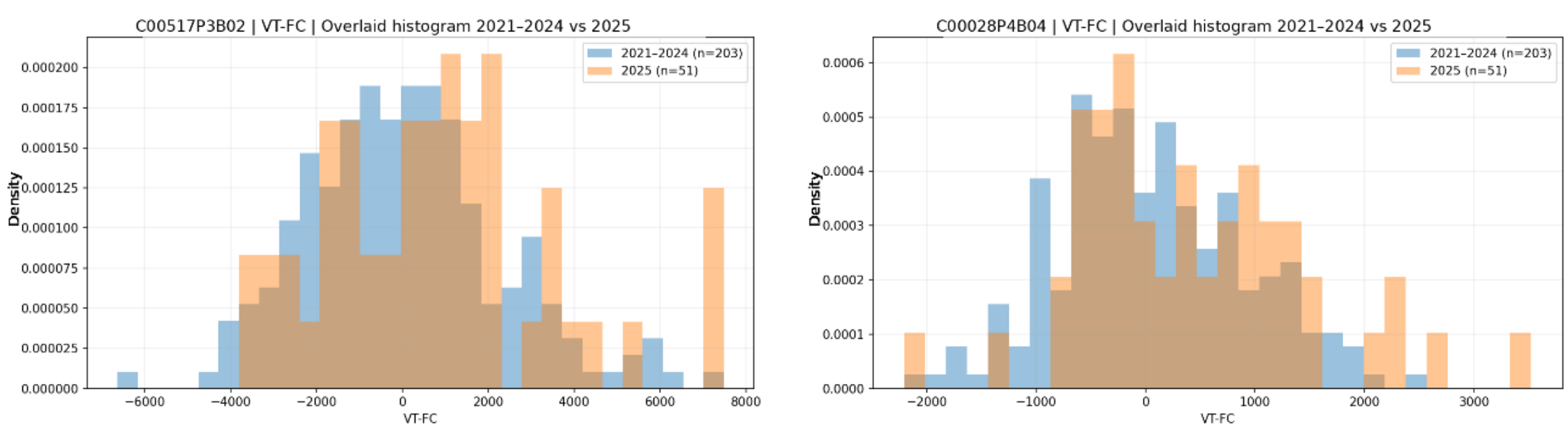

# 4. Conclusions

An adequate historical span is necessary but not sufficient for safety stock estimation. When forecasting processes evolve and learn, the relevant question is not only how many observations to use, but which observations still represent future uncertainty. LOWDII addresses this question by evaluating the distributional contribution of each error to the empirical uncertainty distribution used for safety stock calibration. Our proposed LOWDII method outperforms alternative approaches based solely on temporal proximity to the present, such as exponential smoothing, or on relative distance from the median, such as IQR-based cleansing. In our simulation experiments based on real data, LOWDII achieved the target service level while reducing average stock by between 5.1% and 22.6%, depending on the benchmark method considered. LOWDII offers a way to use longer forecast-error histories without forcing the inventory system to inherit all past distortions. This is particularly relevant in forecasting environments characterized by learning, correction, and evolving commercial behaviour, where the past remains informative but not uniformly representative. The resulting error histories preserve persistent variability while reducing the impact of transient and non-representative extremes. Although the company's current one-year procedure produced statistically comparable results in this validation, a single-year span remains vulnerable to sampling variability and to the specific commercial, promotional, or operational conditions of the immediately preceding year. As shown by the out-of-sample behaviour of the selected SKUs, recent forecast errors may themselves reflect transient distortions rather than persistent uncertainty. LOWDII therefore offers a more robust calibration logic: it exploits a longer historical span while selectively reducing the influence of observations that are distributionally inconsistent with the current uncertainty structure.

## 4.1. Managerial implications

From a managerial perspective, LOWDII is readily applicable in FMCG settings where products have life cycles long enough to provide multi-year forecast-error histories and where the forecasting process is sufficiently stable and disciplined to generate comparable errors over time. LOWDII helps managers translate improvements in the forecasting process into inventory policy earlier. Using it as guidance, firms can convert forecast-error histories into a more representative basis for future uncertainty estimation. This allows firms to capture the economic benefit of reduced forecast uncertainty sooner, while preserving persistent tail behaviour that remains relevant for service protection.

For implementation, firms should maintain systematic traceability of forecast errors and of the commercial or operational conditions under which they were generated. LOWDII can then be incorporated into a periodic review process in which forecast-error histories and safety stock parameters are jointly reassessed. Its diagnostic output helps planners identify observations whose distributional influence is unusually high relative to the SKU-specific error structure and prioritize managerial analysis on exceptional and distortive deviations.

A further managerial implication concerns the behavioural effect of inventory buffers on forecast improvement. If all historical errors are retained indiscriminately and translated into safety stock, the inventory system may absorb even non-representative or correctable forecasting failures. While this reduces the probability of stock-outs in the short term, it can also weaken the organizational incentive to investigate and correct the underlying sources of error. In this sense, excessive protection may teach the system not to break, but also not to learn. LOWDII mitigates this risk by separating errors that should remain part of the uncertainty basis from distortions that should trigger managerial learning rather than permanent inventory protection.

## 4.2. Future research directions

The experimental scope was intentionally limited to isolate LOWDII's incremental value. Future research should evaluate LOWDII with more advanced safety stock models, especially when forecast errors exhibit autocorrelation, heteroscedasticity, heavy tails, or other departures from the assumptions underlying the standard formulation. Although LOWDII helps construct a more representative historical sample, forecast errors should still be independent across periods, have finite variance, and display constant variance over time for the standard safety stock expression to remain valid. When these conditions are violated, cumulative lead-time uncertainty may be misrepresented. Future studies could therefore integrate LOWDII with time-series error modelling, conditional volatility models, non-parametric quantile estimation, or simulation-based safety stock calibration.

Further research should also broaden the inventory-design context in which LOWDII is evaluated, including settings where demand and lead-time uncertainty are jointly relevant, as well as contexts where risk-pooling across products, locations, or echelons materially affects safety stock decisions. Finally, future research should refine LOWDII by exploring adaptive thresholds that account for SKU maturity, error persistence, and the speed at which forecast-error distributions become centred and less dispersed over time.

***Declaration of interests:*** *Luis Fernández-Palacios reports an employment relationship with the company that provided the confidential data used in this study, which is not identified due to commercial confidentiality restrictions. The remaining authors declare that they have no known competing financial interests or personal relationships that could have appeared to influence the work reported in this paper.*

***Declaration of generative AI and AI-assisted technologies in the manuscript preparation process:*** *During the preparation of this work, the authors used Undermind to support the bibliographic search process and identify potentially relevant references. ChatGPT was used to assist in writing Python code for the simulation. All AI-assisted outputs were critically reviewed, verified, and edited by the authors. The authors take full responsibility for the accuracy, integrity, and final content of the manuscript.*

***Data statement:*** *The data used in this study are subject to commercial confidentiality restrictions and are therefore not publicly available. The dataset contains sensitive business information from a private company. The key replicable components of the study are the mathematical derivations, assumptions, and analytical procedures, which are fully described in the manuscript.*

***Funding declaration:*** *No funding was received for the research, authorship, or publication of this article.*

***Author Contributions Statement:*** *Luis Fernández-Palacios contributed to the conception and design of the study, data curation, formal analysis, software development, simulation, interpretation of results, and drafting of the manuscript. Manuel Ceballos contributed to the mathematical and methodological development, validation of the analytical approach, interpretation of results, and critical revision of the manuscript. Yolanda Muñoz-Ocaña contributed to the operations-management framing, managerial interpretation, supervision, and critical revision of the manuscript. All authors approved the final version to be published and agree to be accountable for all aspects of the work.*